\documentclass[twoside]{article}
\usepackage[frenchb]{babel}	% Typo francaise.
\FrenchItemizeSpacingfalse	% Pour avoir des \bullet dans les {itemize}
\usepackage[T1]{fontenc}	% Accents cod\'es dans la fonte.
\usepackage[latin1]{inputenc}	% Accents 8 bits dans le fichier.
\usepackage{vmargin}		% R\'egler la taille de la feuille.
\usepackage{fancyhdr}		% R\'egler le titre courant et le bas de page.
\usepackage{calc}	    	% Faire des calculs sur les longueurs.
\usepackage{ifthen}		    % Faire des tests if/then/else.
\usepackage{pifont}		    % La police \ding.
\usepackage{supertabular}	% Les grands tableaux.
\usepackage{longtable}
\usepackage{multicol}		% Plusieurs colonnes.
\usepackage{wrapfig}		% Dessins dans le texte.
\usepackage{fancybox}		% Boîtes avec une ombre pour les parties.
\usepackage{rotating}		% Pour tourner un texte.
\usepackage{xspace}		    % Ajuster l'espace apr\`es des mots.

\usepackage{amsmath}	% Les symboles les plus fr\'equents.
\usepackage{amssymb}	% Des symboles.
\usepackage{amsfonts}	% Des fontes, eg pour \mathbb.
\usepackage{verbatim}	% Pour les codes sources en informatique.
\usepackage{mathrsfs}	% Des lettres majuscules cursives (\mathscr).
\usepackage{dsfont}

\usepackage{array}
\usepackage{curves}
\usepackage{epic}
\usepackage{eepic}
\usepackage{epsfig}
\usepackage{graphics}
\graphicspath{{PS/}}

\setpapersize{custom}{21cm}{29.7cm}
\setmarginsrb{15mm}{20mm}{15mm}{20mm}{0cm}{10mm}{0cm}{1cm}

\def \no {\noindent}
\def \sm {\smallskip}
\def \d {\textup{d}}
\def \ppcm {\textup{ppcm}}

\newcommand {\td} [1]  {\textup #1}
\newcommand {\ti} [1] {\tilde{#1}}
\newcommand {\Z}  {\mathbb{Z}}
\newcommand {\Q}  {\mathbb{Q}}
\newcommand {\R}  {\mathbb{R}}
\newcommand {\C}  {\mathbb{C}}
\newcommand {\K}  {\mathbb{K}}
\newcommand {\N}  {\mathbb{N}}
\newcommand {\f}  {\frac}
\newcommand {\pa}[1]  {\left( #1\right)}
\newcommand {\pac}[1]  {\left[ #1\right]}
\newcommand {\paa}[1]  {\left\{ #1\right\}}
\newcommand {\abs}[1]  {\left| #1\right|}
\newcommand{\Sum}[2]  {\displaystyle{\sum_{#1}^{#2}}}
\newcommand{\Int}[2]  {\displaystyle{\int_{#1}^{#2}}}
\newcommand{\Prod}[2]  {\displaystyle{\prod_{#1}^{#2}}}
\newcommand{\partieentiere}[1]{\left\lfloor #1 \right\rfloor}
\newcommand{\dt}{\textup{d}t}
\newcommand{\dx}{\textup{d}x}
\newcommand{\dz}{\textup{d}z}
\newcommand{\p}{\mathfrak{p}}

\newtheorem{theorem}{Th\'eor\`eme}

\newtheorem{lem}{Lemme}
\newtheorem{prop}{Proposition}
\newtheorem{cor}{Corollaire}

\begin{document}

\title{Fonction Z\^eta de Hurwitz $p$-adique et irrationalité}
\author{Pierre Bel}
\maketitle

\begin{center}
\begin{minipage}{11cm}

The knowledge on irrationality of $p$-adic zeta values has recently progressed. The irrationality of $\zeta_2(2)$, $\zeta_2(3)$ and of a few other $p$-adic series of Dirichlet  was obtained by F. Calegari(cf. \cite{Ca}). F. Beukers gave a more elementary proof of this result(cf. \cite{Be}). In parallel, T. Rivoal has just obtained, in the complex case, some Pad\'e approximants of Lerch functions (cf. \cite {Ri2}). It is this work which, transposed to $\C_p$, enables us to obtain results of irrationality and linear independence.

\end{minipage}
\end{center}
\section{Introduction}

\subsection{Pr\'eliminaires}
\ \ \ Soit $p$ un nombre premier. On note $v_p$ la valuation $p$- adique sur $\Q$ et $|.|_p=p^{-v_p}$ la valeur absolue $p$-adique. On pose $q_p=p$ si $p\not=2$ et $q_2=4$. Pour $x\in\Z_p^*$, on d\'esigne par $\omega(x)$ l'unique racine de l'unit\'e, d'ordre $p-1$ si $p \not=2$, et d'ordre $2$ si $p=2$, telle que $\vert x-\omega(x)\vert_p \le q_p^{-1}$. On \'etend la d\'efinition de $\omega$ \`a $\Q_p^*$ en posant $\omega (x)=p^{v_p(x)}\omega(p^{-v_p(x)}x)$, et on pose $<x>={x \over\omega(x)}$ (donc $<px>=<x>$ pour tout $x\in\Q_p^*$).

\ \ \ On note $\log_p$ la fonction d\'efinie par $$\log_p(1+x)=\sum_ {k=1}^{+\infty}(-1)^{k+1}{x^k\over k}$$ pour $x\in\C_p$ tel que $ \vert x\vert_p<1$.

\ \ \ 0n note $\zeta(s,x)$ la fonction z\^eta de Hurwitz d\'efinie par
$$
\zeta(s,x)=\sum_{k=0}^{+\infty}{1\over(n+x)^s}
$$
pour $(s,x)\in\C\times\R$, avec $\Re(s)>1$ et $x>0$. Pour $x$ fix\'e, cette fonction admet un prolongement en une fonction holomorphe sur $ \C\backslash\{1\}$, dont $1$ est un p\^ole d'ordre $1$ et de r\'esidu $1$.\par\ \ \ La formule d'Euler-MacLaurin conduit ais\'ement au d\'eveloppement asymptotique suivant, pour $x\rightarrow+\infty$ ~:
\begin{eqnarray}\label{zetaasympt1}
\zeta(s,x)={x^{1-s}\over s-1}-\sum_{j=1}^k{-s\choose{j-1}}{B_j\over j} x^{1-s-j}+O(x^{-s-k})
\end{eqnarray}

\no o\`u les $B_j$ sont les nombres de Bernoulli, et pour tout $k\ge1$ le symbole $O$ est uniforme en $s$ pour $s$ born\'e. Par passage \`a la limite sur $s$, on en d\'eduit que la valeur en 1 de la fonction holomorphe $s\longmapsto\zeta(s,x)-{1\over s-1}$ v\'erifie~:
\begin{eqnarray}\label{zetaasympt2}
\left(\zeta(s,x)-{1\over s-1}\right)_{s=1}=-\ln x+\sum_{j=1}^k{(-1) ^jB_j\over j}x^{-j}+O(x^{-k-1}).\end{eqnarray}

\no
La fonction z\^eta $p$-adique de Hurwitz peut \^etre exprim\'ee par son d\'eveloppement en s\'erie de Laurent~:

\begin{eqnarray}\label{zetapasympt1}
\zeta_p(s,x)={<x>^{1-s}\over s-1}-<x>^{1-s}\sum_{j=1}^{+\infty}{-s \choose{j-1}}{B_j\over j}x^{-j}.
\end{eqnarray}

\no
Ce d\'eveloppement est $p$-adiquement convergent pour $\vert x\vert_p >1$ car le nombre $$\binom {-s}{j-1}=(-1)^{j-1}\binom{s+j-2}{j-1}$$ est entier. On a aussi~:

\begin{eqnarray}\label{zetapasympt2}
\lim_{s\rightarrow1}(\zeta_p(s,x)-{1\over s-1})=-\log_p <x>+\sum_{j=1}^ {+\infty}{(-1)^jB_j\over j}x^{-j}.\end{eqnarray}

\no
On pourra se r\'ef\'erer au livre de H. Cohen   (\cite{Coh}) pour une vision exhaustive de ces diff\'erents r\'esultats.

\no

  \vspace{0.5cm}
  La connaissance sur l'irrationalit\'e de valeurs des fonctions z \^eta $p$-adiques a progress\'e r\'ecemment. L'irrationalit\'e de $ \zeta_2(2)$, $\zeta_2(3)$ et de quelques autres s\'eries de Dirichlet $p$-adiques a \'et\'e obtenue
par F. Calegari(cf. \cite{Ca}). F. Beukers en a donn\'e une interpr\'etation plus \'el\'ementaire(cf. \cite{Be}). Parall\`element, T. Rivoal vient d'obtenir, dans le cas complexe,
certains approximants de Pad\'e de fonctions de Lerch (cf. \cite {Ri2}) et d'\'etudier leur propri\'et\'es diophantiennes. C'est ce travail qui, transpos\'e \`a $\C_p$, nous permet d'obtenir des r\'esultats d'irrationalit\'e et d'ind\'ependance lin\'eaire, gr\^ace \`a un crit\`ere comparable \`a celui de Nesterenko, mais dans lequel nous utilisons des formes lin\'eaires suppos\'ees a priori ind\'ependantes. Un point crucial de notre travail sera d'ailleurs de v\'erifier que cette condition est bien satisfaite dans l'application que nous en ferons (N.D.L.R. lemme du d\'eterminant).
\subsection{R\'esultats}

Soient un entier $e\geq 2$, $v=\ppcm(e,p-1)$, $\xi$ une racine primitive $e$-\`eme de l'unit\'e et $\chi$  une racine primitive $v$-\`eme de l'unit\'e.

\no
Pour un nombre $p$-adique $x$, tel que $|x|_p\geq p$ et $s$ un entier strictement positif, on pose

\begin{eqnarray}
\tilde T_p(s,x)&=\displaystyle{\sum_{j=0}^{e-1}\xi^{-j}\zeta_p\pa{s, \frac {x+j} e}} &\textup{ si }s>1 
\text{ et }\abs{x}\geq q_p\\
&= \displaystyle{\sum_{j=0}^{e-1}((-1)^{s-1}\xi)^{-j}\zeta_2\pa{s, \frac {x+j} e}} &\textup{ si }s>1,\,\ p=2 
\text{ et }\abs{x}_2= 2
\end{eqnarray}
et

\begin{eqnarray}
\tilde  T_p(1,x)&= \displaystyle{\lim_{s\rightarrow1}\tilde T_p(s,x)}=\displaystyle{\lim_{s\rightarrow1}\sum_{j=0}^{e-1} \xi^{-j}{\zeta_p(s,\f {x+j} e)} }.
\end{eqnarray}

\begin{theorem}

Soit $x=\f{a}b$ un rationnel, tel que $|x|_p\geq p$, et soit $A$ un entier sup\'erieur ou \'egal \`a 2.

\no
Alors la dimension $\tau$ de l' espace vectoriel engendr\'e sur $\Q (\chi)$ par la famille
$\paa{1,\left(\tilde T_p(s,x)\right)_{s\in[1,A]}}$ v\'erifie

\medskip

$$\tau\geq\frac{[\Q_p(\xi):\Q_p]}{\varphi(v)}\frac{A\ln\abs{x}_p}{\ln b+\Sum{q| b}{}\f{\ln q}{q-1}+
A+(A-1)\ln2}.$$

\end{theorem}

\begin{theorem}
Pour tout entier $A$ sup\'erieur ou \'egal \`a $2$, il existe une borne $M_A$ explicite tel que si le nombre premier $p$ est plus grand que $M_A$ alors

\begin{center}{l'ensemble $\paa{\zeta_p\pa{s, \f 1 p}-\zeta_p\pa{s, \f {p+2}{ 2p}}}_{s\in[1,A]}$ contient au moins $A-1$ nombres irrationnels.}
\end{center}

\end{theorem}

Le point crucial de la d\'emonstration des r\'esultats est le calcul du d\'eterminant de la partie 5 qui permet d'appliquer le crit\`ere d'ind\'ependance lin\'eaire suivant.

\section{Crit\`ere d'ind\'ependance lin\'eaire}

On rappelle la formule du produit pour un corps de nombres $\K$. Pour $v$ une place de $\K$, on note $K_v$ et $Q_v$ les compl\'et\'es de $\K $ et $\Q$ pour cette place et $\eta_v=[\K_v:\Q_v]$.

\no
Si $\alpha\in\K^*$, alors on a

$$0=\sum_{v\text{ place de } \K}\eta_v\ln\abs{\alpha}_v.$$

\no
De plus

$$\sum_{v\text{ place de } \K\text{ infinie }}\eta_v=[\K:\Q]$$

\no
Si $\alpha$ est un \'el\'ement non nul
de ${\mathcal O}(\K)$, comme $\vert\alpha\vert_v\le 1$ pour toute place finie $v$ de $\K$, si $\p$ est une place finie, on  a~:
$$\eta_\p\ln\vert\alpha\vert_\p+\sum_{v\ \rm infinie}\eta_v\ln\vert \alpha\vert_v\ge0.$$
\ \ \ Soit $m$ un nombre entier positif. Pour $L= (\ell_1,...\ell_m) \in\K^m$, et $\theta=(\theta_1,...,\theta_m)\in \C_p^m$, on note $(L, \theta)=\ell_1\theta_1+...+\ell_m\theta_m$.\par
\no
Si $v$ est une place de $\K$, on note $\Vert L\Vert_v=\max_{1 \le j \le m}\vert\ell_j\vert_v$.

\vspace{0.3cm}

\no

\begin{lem} Soit $p$ un nombre premier. Soit $\K$ un corps de
nombres sur
$\Q$. On consid\`ere $\K$ comme plong\'e dans $\C_p$,  dans lequel on  note
$\K_p=\Q_p(\K)$ son adh\'erence. Soit $\theta= (\theta_1,\cdots, \theta_m)$ un
\'el\'ement non nul de $\K_p$. On
suppose qu'il existe $m $ suites $(L_n^{(i)})$, o\`u $n\in\N$ et
$\ 1\le i\le m$, d'\'el\'ements de $(\mathcal{O}(\K))^m$ telles que pour chaque $n$ les $L_n^ {(i)}$, pour $1\le i\le m$,
soient lin\'eairement ind\'ependants sur $\K$, et des nombres r\'eels strictement positifs $c$ et $\rho$
satisfaisant pour chaque $i$~ les conditions~:
$$\limsup_n\f1n\ln{\Vert L_n^{(i)}\Vert_v}\leq c$$
\no
pour les places infinies $v$, et
$$\limsup_n\f1n\ln\abs{(L^{(i)}_n,\theta)}_p\leq-\rho.$$
\no
Alors la dimension $\tau$ du
$\K$-sous-espace vectoriel de $\K_p$ engendr\'e par les $\theta_j$ pour $1\le j\le m$ v\'erifie

$$\tau\geq \frac{ \rho\, [\K_p:\Q_p]}{c\,[\K:\Q]}.$$

\end{lem}

\no
\textbf{D\'emonstration}

Effectuons tout d'abord quelques r\'eductions. En renum\'erotant les variables $(\theta_i)_{i\in[1,m]}$, on
peut supposer $\theta_1\neq 0$. De plus, en remplaçant les variables $ (\theta_j)_{j\in[1,m]}$, par $\left(\f{\theta_j}{\theta_1}\right)_{j \in[1,m]}$, les hypoth\`eses \'etant encore v\'erifi\'ees, on peut  supposer $\theta_1=1$.

\sm

Si $\tau$ est la dimension  du $\K$-espace vectoriel engendr\'e par les $\theta_j$, alors il existe $m-\tau$ \'el\'ements $(A^{(i)})_{i\in [\tau+1,m]}$ de $(\mathcal{O}(\K))^m$, lin\'eairement ind\'ependants  sur $\K$, tels que  $
(A^{(i)},\theta)=0$ pour tout $i\in[\tau+1,m]$.

\sm

On peut en faisant des permutations entre les $L_n^{(i)}$ \`a chaque  rang
$n$ se ramener au cas, o\`u pour tout $n\in\N$, la famille $(L_n^ {(1)},\cdots,L_n^{(\tau)},A^{(\tau+1)},\cdots,A^{(m)})$ est libre.

Soit $M_n$ la matrice dont les lignes sont form\'ees des vecteurs
$(L_n^{(1)},\cdots,L_n^{(\tau)},A^{(\tau+1)},\cdots,A^{(m)})$, i.e., en posant $L_n^{(i)}= (\ell_{n,1}^{(i)},...\ell_{m,1}^{(i)},)$ et $A^ {(i)}=(a^{(i)}_1,\cdots,a^{(i)}_m)$,

$$M_n=
\left(
\begin{array}{cccc}
\ell_{n,1}^{(1)}& \ell_{n,2}^{(1)}& \cdots & \ell_{n,m}^{(1)}\\
\cdots & \cdots& \cdots& \cdots\\
\ell_{n,1}^{(\tau)}& \ell_{n,2}^{(\tau)}& \cdots & \ell_{n,m}^ {(\tau)} \\
a^{(\tau+1)}_1 & a^{(\tau+1)}_2 & \cdots &  a^{(\tau+1)}_m \\
\cdots & \cdots& \cdots& \cdots\\
a^{(m)}_1 & a^{(m)}_2 & \cdots & a^{(m)}_m \\

\end{array}
\right)$$

\no
Comme la matrice est non singuli\`ere, on a

$$\Lambda_n=\det(M_n)\neq0$$

\no
Comme $\Lambda_n$ appartient \`a $\mathcal{O}(\K)$, on en d\'eduit que

\begin{eqnarray}\label{ineglambda2}
0\leq \eta_p\ln|\Lambda_n|_p+\sum_{v\text{ infinie }} \eta_{v}\ln| \Lambda_n|_v.
\end{eqnarray}

\no
Le d\'eveloppement du d\'eterminant, nous permet
d'obtenir pour les places infinies:

\begin{eqnarray}\label{majinflambda}\limsup_n\f{\ln{\abs{\Lambda_n} _v}}n\leq \tau \,c.
\end{eqnarray}

\sm

\no
Pour le calcul du d\'eterminant, on peut aussi ajouter \`a une colonne, une combinaison
lin\'eaire des autres colonnes. En ajoutant \`a la premi\`ere, les colonnes suivantes
respectivement multipli\'ees par $\theta_j$, on obtient:

$$\Lambda_n=
\left|
\begin{array}{cccc}
(L_n^{(1)},\theta)&\ell_{n,2}^{(1)} & \cdots & \ell_{n,m}^{(1)}\\
\cdots & \cdots& \cdots& \cdots\\
(L_n^{(\tau)},\theta)& \ell_{n,2}^{(\tau)}& \cdots & \ell_{n,m}^ {(\tau)} \\
0 & a^{(\tau+1)}_2 & \cdots & a^{(\tau+1)}_m \\
\cdots & \cdots& \cdots& \cdots\\
0 & a^{(m)}_2& \cdots & a^{(m)}_m \\

\end{array}
\right|$$

\sm

\no
Le d\'eveloppement du d\'eterminant sous cette forme nous permet
d'obtenir:

\begin{eqnarray}\label{majplambda}
\limsup_{n}\f{\ln{\abs{\Lambda_n}_p}}n\leq -\rho
\end{eqnarray}

\no
En divisant l'in\'equation (\ref{ineglambda2}) par $n$ et utilisant (\ref{majinflambda}) et (\ref{majplambda}), on en
d\'eduit:

$$0\leq -\rho\,\eta_p +\tau c\sum_{v\text{ infinie }}\eta_v$$

\no
Comme  $\displaystyle{\sum_{v\text{ infinie }}\eta_v=[\K:\Q]}$, le r\'esultat est donc d\'emontr\'e.

\sm

\section{Approximants de Pad\'e simultan\'es de fonctions Z\^eta de  Hurwitz}
Dans cette partie, comme dans la suite, $\xi$ est une racine primitive $e$-\`eme de l'unit\'e avec $e$ entier, $e\ge2$.

\no
On d\'efinit, pour $x$ un nombre complexe diff\'erent d'un entier n\'egatif, $s$ un
entier strictement positif et $z$ un nombre complexe tel que
$\vert z\vert\le1$ et $z\not=1$, $\phi$ la fonction de Lerch~:

$$\phi_s(x,z)=\sum_{k=0}^{+\infty}\f{z^k}{(k+x)^s}$$
\no
On remarque pour $\Re(x)>0$ et pour $m>0$ entier, l'expression

$$\phi_m(x,z)=\sum_{k=0}^{+\infty}\f{z^k}{(k+x)^m}=\sum_{k=0}^{+ \infty}\frac{(-1)^{m-1}}{(m-1)!}z^k \int_0^1 t^{x-1+k}(\ln t)^{m-1}\dt= \frac{(-1)^{m-1}}{(m-1)!} \int_0^1\frac {t^{x-1}(\ln t)^{m-1}}{1-zt} \dt$$
\no
qui montre la convergence de la s\'erie
$\sum_{k=0}^{+\infty}\f{z^k}{(k+x)^m}$ pour  
$\vert z\vert\le1$ et $z\not=1$, et 
permet de prolonger la fonction $\phi_m(x,z)$ en une fonction holomorphe en $z$ sur $\C\setminus[1,+\infty[$. Par translation enti\`ere, il en est
finalement ainsi pour tout nombre complexe $x$ tel que $-x\not\in\N$. 

On suppose que $A$ est un entier sup\'erieur ou \'egal \`a 2, consid\'er\'e comme fix\'e. Le nombre $n$ est un entier
positif v\'erifiant $A\,n\geq n+3$.

On rappelle que le symbole de Pochamer est not\'e
$$(t)_m=\prod_{0\le j<m}(t+j)$$
pour $t\in\C$ et $m\in\Z$.

Posons pour $q\in[0,A]$ et un nombre $x$ tel que $x\notin \Z^-$

$$R^{(q)}_n(k)=n!^{A-1}\f{(k)_{n+1}}{(k+x)^A_{n}(x+k+n)^q}$$

\no
et
$$S^{(q)}_{n}(x,z)=\sum_{k=0}^{\infty}{R^{(q)}_n(k)\,z^{-k}}.$$

La fraction rationnelle $R^{(q)}_n(k)$ est de degr\'e $n+1-A\,n-q$ par rapport \`a $k$
donc de degr\'e inf\'erieur ou \'egal \`a $-2$, vu les hypoth\`eses. La fonction
$S^{(q)}_{n}(x,z)$ est donc d\'efinie pour tout complexe $z$  de module sup\'erieur
ou \'egal \`a $1$. La s\'erie $S^{(q)}_{n}(x,z)$ converge normalement sur l'ensemble des complexes $x$ de partie r\'eelle plus grande que $1$ et des complexes $z$ de module plus grand que $1$.
\sm

\sm

\begin{prop} \label{approxpade}
Il existe  $A+1$ polyn\^omes $(P^{(q)}_{s}(x,z))_{s\in[0,A]}$  \`a coefficients
rationnels de degr\'e en $x$ au plus $n+1$, de degr\'e en $z$ au plus $n$, et le degr\'e en $z$
de $P^{(q)}_s(x,z)$ est  m\^eme au plus $n-1$ si $s>q$, tels que pour tout $z$ avec $\vert z\vert\ge1$ et $z\not=1$
et tout $x\not\in-\N$, on ait~:

\begin{eqnarray}\label{approxpade1}
S^{(q)}_{n}(x,z)=P^{(q)}_0(x,z)+\sum_{s=1}^A P^{(q)}_{s}(x,z)\phi_s \left(x,\f1z\right).
\end{eqnarray}

De plus, on a, quand $\Re(x)\to+\infty$
\sm

$$S^{(q)}_{n}(x,z)=o(x^{-A\,n+n+3-q}).$$

\end{prop}
\sm

\textbf{D\'emonstration}

La d\'ecomposition en \'el\'ements simples de $R^{(q)}_n(k)$ nous donne
$$R^{(q)}_n(k)=\sum_{s=1}^A\sum_{j=0}^n \f{r^{(q)}_{j,s}(x)}{(k+x+j) ^s}$$

o\`u

\sm

$$r^{(q)}_{j,s}(x)=\left\{
\begin{array}{cl}
\f1{(A-s)!}\left(\f{\td{d}}{\td{d}k}\right)^{A-s}
\left[R^{(q)}_n(k)(x+k+j)^A\right]_{|k=-j-x} &\textup{si $j\in [0,n-1]$ et $s\in[1,A]$}\\[4mm]
\f1{(q-s)!}\left(\f{\td{d}}{\td{d}k}\right)^{q-s}
\left[R^{(q)}_n(k)(x+k+n)^q\right]_{|k=-n-x} & \textup{si $j=n$ et $s\in[1,q]$}\\[4mm]
0& \textup{si $j=n$ et $s\in[q+1,A]$}\\
\end{array}\right..$$

\sm
Remarquons tout de suite que, pour $q>0$,

\begin{eqnarray}\label{rqq}
r^{(q)}_{n,q}(x)=\left[R^{(q)}_n(k)(x+k+n)^q\right]_{k=-n-x}=n!^ {A-1}\f{(-n-x)_{n+1}}{(-n)^A_{n}}\neq0
.\end{eqnarray}

\sm

Par le changement de variable $l=-k-x$, on obtient

$$r^{(q)}_{j,s}(x)=\left\{
\begin{array}{cl}
\f{(-1)^{A-s}}{(A-s)!}\left(\f{\td{d}}{\td{d}l}\right)^{A-s}
\left[R^{(q)}_n(-l-x)(j-l)^A\right]_{|l=j} &\textup{si $j\in[0,n-1] $ et $s\in[1,A]$}\\[4mm]
\f{(-1)^{q-s}}{(q-s)!}\left(\f{\td{d}}{\td{d}l}\right)^{q-s}
\left[R^{(q)}_n(-l-x)(j-l)^q\right]_{|l=n} & \textup{si $j=n$ et $s \in[1,q]$}\\[4mm]
0& \textup{si $j=n$ et $s\in[q+1,A]$}\\
\end{array}\right..$$

On en d\'eduit

\begin{eqnarray}\label{rjs}
r^{(q)}_{j,s}(x)=\left\{
\begin{array}{cl}
\f{(-1)^{A-s}}{(A-s)!}\left(\f{\td{d}}{\td{d}l}\right)^{A-s}
\left[n!^{A-1}\f{(-l-x)_{n+1}}{(-l)^A_{n}(n-l)^q}(j-l)^A\right]_{| l=j} &\textup{si $j\in[0,n-1]$ et $s\in[1,A]$}\\[4mm]
\f{(-1)^{q-s}}{(q-s)!}\left(\f{\td{d}}{\td{d}l}\right)^{q-s}
\left[n!^{A-1}\f{(-l-x)_{n+1}}{(-l)^A_{n}}\right]_{|l=n} & \textup{si $j=n$ et $s\in[1,q]$}\\[4mm]
0& \textup{si $j=n$ et $s\in[q+1,A]$}\\
\end{array}\right..
\end{eqnarray}

\no
Les fonctions $r^{(q)}_{j,s}(x)$ sont donc des polyn\^omes en $x$ de degr\'e au
plus $n+1$.

$$S^{(q)}_{n}(x,z)=\sum_{k=0}^{\infty}
{\sum_{s=1}^A\sum_{j=0}^n \f{r^{(q)}_{j,s}(x)}{(k+x+j)^s}z^{-k}}.$$

\no
Il en r\'esulte que

$$
\begin {array}{lcll}
  S^{(q)}_{n}(x,z) &= & \Sum{s=1}A\Sum{j=0}n\Sum{k=0}{+\infty}
\f{r^{(q)}_{j,s}(x)}{(k+x+j)^s}z^{-k}& \\[4mm]
   &= &\Sum{s=1}A\Sum{j=0}n r^{(q)}_{j,s}(x)z^j\Sum{k=0}{+\infty}
\f{z^{-k-j}}{(k+x+j)^s}\\[4mm]
   & =&\Sum{s=1}A\Sum{j=0}n r^{(q)}_{j,s}(x)z^j\left[\phi_s\left(x,\f1z\right)- \Sum{k=0}{j-1}
\f{z^{-k}}{(k+x)^s}\right]\\[4mm]
&=&\Sum{s=1}A\phi_s\left(x,\f1z\right)\Sum{j=0}n r^{(q)}_{j,s}(x)z^j-\Sum{s=1}A \Sum{j=0}n
r^{(q)}_{j,s}(x)z^j\Sum{k=0}{j-1}
\f{z^{-k}}{(k+x)^s}.\\[10mm]
\end{array}.
$$

\no
On a donc
$$
\begin{array}{lcll}
S^{(q)}_{n}(x,z) =P^{(q)}_0(x,z)+\Sum{s=1}A P^{(q)}_s(x,z)\phi_s\left (x,\f1z\right)\\
\end{array},
$$

\no
\textup{ o\`u l'on a pos\'e}
\sm

\begin{center}{$P^{(q)}_0(x,z)=-\Sum{s=1}A\Sum{j=0}n
r^{(q)}_{j,s}(x)z^j\Sum{k=0}{j-1}\f{z^{-k}}{(k+x)^s}$},
\end{center}

\no
\textup{et, pour tout $ s \in[1,A]$}
\begin{eqnarray}\label{Pqs}
P^{(q)}_s(x,z)=\Sum{j=0}n r^{(q)}_{j,s}(x)z^j.
\end{eqnarray}

\sm

Les \'egalit\'es (\ref{rjs}) montrent imm\'ediatement que pour
$s\geq1$, $P^{(q)}_{s}(x,z)$ est un polyn\^ome \`a coefficients rationnels de degr\'e
en $x$ au plus $n+1$ et de degr\'e en $z$ au plus $n$. On voit directement que le degr\'e en $z$ de $P^{(q)}_s$ est au plus $n-1$, si $s>q$.

Pour $P^{(q)}_0(x,z)$, on voit directement que c'est un polyn\^ome
en $z $ de degr\'e
au plus $n$. De plus, on remarque que

\begin{center}{$\Sum{k=0}{j-1}\f
{z^{j-k}}{(k+x)^s}=\f{(-1)^{s-1}}{(s-1)!}\left(\f{\td{d}}{\td{d}l} \right)^{s-1}\pac{\Sum{k=1}{j}\f
{z^{k}}{(l-k+x)}}_{|l=j}$}
\end{center}
\sm

\no
Il en r\'esulte que pour $j\in[1,n-1]$

\begin{eqnarray*}\begin{array}{ccl}
\Sum{s=1}A r^{(q)}_{j,s}(x)\Sum{k=0}{j-1}\f{z^{j-k}}{(k+x)^s}
&=&\Sum{s=1}A\f{(-1)^{s-1}}{(s-1)!}\left(\f{\td{d}}{\td{d}l}\right)^ {s-1}\pac{\Sum{k=1}{j}\f
{z^{k}}{l-k+x}}_{|l=j}\f{(-1)^{A-s}}{(A-s)!}\left(\f{\td{d}}{\td{d}l} \right)^{A-s}
\left[R^{(q)}_n(-l-x)(j-l)^A\right]_{|l=j}\\[6mm]
&=&\f{(-1)^{A-1}}{(A-1)!}\Sum{s=1}A\binom{A-1}{s-1}\left(\f{\td{d}}{\td {d}l}\right)^{s-1}\pac{\Sum{k=1}{j}\f
{z^{k}}{l-k+x}}_{|l=j}\left(\f{\td{d}}{\td{d}l}\right)^{A-s}
\left[R^{(q)}_n(-l-x)(j-l)^A\right]_{|l=j}\\[6mm]
\end{array}
\end{eqnarray*}

\begin{eqnarray}\begin{array}{ccl}
&=&\f{(-1)^{A-1}}{(A-1)!}\left(\f{\td{d}}{\td{d}l}\right)^{A-1}
\left[R^{(q)}_n(-l-x)(j-l)^A\Sum{k=1}{j}\f {z^{k}}{l-k+x}\right]_{| l=j}\label{rfrac1}
\end{array}
\end{eqnarray}

\no
On a

\begin{eqnarray}\label{rfrac2}
R^{(q)}_n(-l-x)(j-l)^A\Sum{k=1}{j}\f
{z^{k}}{l-k+x}=n!^{A-1}\f{(-l-x)_{n+1}}{(-l)^A_{n}(-l+n)^q}(j-l)^A
\Sum{k=1}{j}\f
{z^{k}}{l-k+x}
\end{eqnarray}

\no
Comme les p\^oles simples  en $x$ de $\Sum{k=1}{j}\f {z^{k}}{l-k+x}$
sont des z\'eros de $(-l-x)_{n+1}$, 
$R^{(q)}_n(-l-x)(j-l)^A\Sum{k=1}{j}\f {z^{k}}{l-k+x}$ est un polynôme en $x$
de degr\'e au plus $n$. On justifie de mani\`ere similaire le cas $j=n $ et il en
r\'esulte que $P^{(q)}_0(x,z)$ est un polyn\^ome de degr\'e au plus $n $ par rapport \`a
$x$. La premi\`ere partie de la proposition est donc d\'emontr\'ee.

\medskip

Pour le dernier point, on a la majoration pour $\Re(x)>0$:

$\left|x^{An-n-3+q}S^{(q)}_{n}(x,z)\right|\leq n!^{A-1}\Sum{k=0}{+ \infty}
\f{(k)_{n+1}|x|^{An-n-3+q}}{|k+x|^{An+q}}$

\begin{center}{$\leq n!^{A-1}\Sum{k=0}{+\infty}
\f{(k)_{n+1}|x|^{An-n-3+q}}{|x+k|^{An-n-3+q}|k+x|^{n+3}}\leq n!^{A-1} \Sum{k=0}{+\infty}
\f{(k)_{n+1}}{|k+x|^{n+3}}$}
\end{center}

La convergence normale de la derni\`ere s\'erie sur l'ensemble des complexes $x$ tels que $\Re(x)>1$ permet de passer \`a la limite sous le signe somme
et on conclut que

$$\lim_{\Re(x)\rightarrow+\infty}{\abs{x^{An-n-3+q}S^{(q)}_{n}(x,z)}} =0.$$

La proposition est donc d\'emontr\'ee.
\sm

\begin{cor}\label{corcombSxxi} On a

\begin{center}{$S^{(q)}_{n}(x,\xi)=P^{(q)}_0(x,\xi)+\Sum{s=1}{A}P^ {(q)}_s(x,\xi)\phi_s(x,\xi^{-1})$}
\end{center}

\sm

et, lorsque $\Re(x)\rightarrow+\infty$,
\begin{center}{$S^{(q)}_{n}(x,\xi)=o(x^{-An+n+3})$}
	\end{center}
\end{cor}

\begin{lem} On a:

$$\phi_s(x,1)=\zeta\pa{s,x}$$

$$\phi_s(x,\xi^{-1})=\f1{e^s}\Sum{j=0}{e-1}{\xi^{-j}\zeta\pa{s,\f{x+j}e}}$$

\end{lem}
La preuve est \'evidente.
\vspace{.5cm}

\sm

\section{Propri\'et\'es arithm\'etiques des polyn\^omes $P_s^{(q)}(x,z)$}
\no
On pose $d_n=\ppcm(1,\cdots,n)$. On sait par le th\'eor\`eme des nombres premiers que

$$\ln d_n\sim n.$$

\no
On pose pour tout entier $b$ non nul et pour entier positif $n$
$$\mu_n(b)=b^n\prod_{q|b}q^{\partieentiere{ \f n
{q-1}}}.$$
(o\`u $q$ d\'esigne un nombre premier).

\begin{lem} \label{aritmu}
Si $x$ est un nombre rationnel $\f a b$ ($b>0$) et $k$ un entier appartenant \`a l'intervalle $[0,n]$, alors
les nombres $\displaystyle\f{(x)_n}{n!}\mu_n(b)\text{ et }\f{(x)_{n +1}}{n! (x+k)}\mu_n(b)d_n\text{ sont des entiers }$ et on a

\begin{eqnarray}\label{limmu}{
\lim_{n\rightarrow+\infty}\f1n\ln{(\mu_{n}(b))}=\ln b+\sum_{q|b}\f {\ln q}{q-1}}
.\end{eqnarray}

\end{lem}
\sm

\no
\textbf{D\'emonstration}

\no
On a

$$\f{(x)_n}{n!}\mu_n(b)=\f{\Prod{i=0}{n-1}{(bi+a)}}{n!}\prod_{q|b}q^ {\partieentiere{ \f n
{q-1}}}.$$
\no
Montrons que la valuation $q$-adique de ce nombre rationnel
est positive ou nulle pour tout nombre premier $q$.
\begin{itemize}
\item Si $q$ divise $b$, alors la valuation $q$-adique de $n!$ \'etant au plus $\partieentiere{ \f n{q-1}}$, on en d\'eduit que la valuation $q$-adique est positive ou nulle.

\item Si $q$ ne divise pas $b$, alors la valuation $q$-adique de $ \Prod{i=0}{n-1}{(bi+a)}$ est \'egale \`a celle de $\Prod{i=0}{n-1}{(i+ \f a b)}$. Dans l'intervalle $[0,n-1]$, pour un entier positif $j$,  il y a  au moins $\partieentiere{ \f n{q^j}}$ entiers congrus \`a $- \f a b$ modulo $q^j\Z_q$. La valuation $q$-adique de $\Prod{i=0}{n-1} {(i+\f a b)}$ est donc au moins $\Sum{j=1}\infty{}\partieentiere{ \f n
{q^j}}$ qui est exactement la valuation $q$-adique de $n!$. La valuation $q$-adique est donc positive ou nulle.
\end{itemize}

\vspace{.3cm}
\no
Le nombre $\displaystyle{\f{(x)_n}{n!}\mu_n(b)}$ est donc bien un  entier.

\vspace{.5cm}
\no
On a

$$\f{(x)_{n+1}}{n!(x+k)}\mu_n(b) d_n=\f{\displaystyle{\prod_{0\leq i \leq n, i\neq k}{(bi+a)}}}{n!}\pa{\prod_{q|b}q^{\partieentiere{ \f n
{q-1}}}}d_n.$$

\vspace{.5cm}
\no
Pour cela, montrons que pour tout nombre premier $q$, la valuation $q $-adique de ce nombre rationnel est positive ou nulle.

\no
Si $q$ divise $b$, ceci est \'evident puisque $v_q(n!)<\f n{q-1}$.\par
On suppose donc que $q$ ne divise pas $b$.
Pour tout entier $j$ compris entre $1$ et $J=\partieentiere{\frac{\ln n}{\ln q}}$, on d\'esigne par $\nu_j$ le nombre d'entiers $i$ v\'erifiant  $0\leq i \leq n$, $i\neq k$ et $i\equiv -\frac a b$ mod $q^j$. Le nombre

$$Y=\prod_{0\leq i\leq n, i\neq k}^{n}{(bi+a)}$$

\no
est de valuation $q$-adique
$$v_q(Y)\geq \sum_{j=1}^{J-1}j(\nu_{j}-\nu_{j+1})+J\nu_J=\sum_{j=1}^ {J}\nu_j.$$

\no
Pour chaque $j$ compris entre $1$ et $J$, et pour chaque entier tel que $0\leq K\leq \frac n {q^j}-1$, il y a un entier $i$ appartenant \`a l'intervalle $[Kq^j,(K+1)q^j[$ tel que $i\equiv -\frac a b$ mod $q^j$. Le nombre de ces intervalles disjoints est $\partieentiere {\frac n {q^j}}$, par suite $\nu_j\geq\partieentiere{\frac n {q^j}}-1$. On a donc

$$v_q(Y)\geq \sum_{j=1}^{J}\partieentiere{\frac n {q^j}}-J.$$

\no
Or $v_q(n!)=\sum_{j=1}^{J}\partieentiere{\frac n {q^j}}$ et $v_q(d_n)=J$, on en d\'eduit

$$v_q(Y)-v_q(n!)+v_q(d_n)\geq 0.$$

\no
On conclut que le nombre $\f{(x)_{n+1}}{n!(x+k)}\mu_n(b) d_n$ est de valuation $q$-adique positive ou nulle.
Le nombre $\f{(x)_{n+1}}{n! (x+k)}\mu_n(b)d_n$ est donc bien un entier.

\vspace{.5cm}
\no
Pour la limite (\ref{limmu}), le calcul est direct.
\vspace{.5cm}

\begin{prop}\label{coefentier}
Pour tout nombre premier $p$, et tout $s\in[1,A]$, on a
$$p^{\partieentiere{\f{n}{p-1}}}d_n^{A-s}\, P^{(q)}_s(x,\xi)\in\Z_p [\xi][x] $$
\no
et
$$p^{\partieentiere{\f{n}{p-1}}}d_n^{A-1}\, P^{(q)}_0(x,\xi)\in\Z_p [\xi][x] .$$
\end{prop}

\no De plus, pour un nombre rationnel $\f a b$, avec $(a,b)=1$, pour tout $s\in[1,A]$,
on a 

$$ b\,d_n^{A-s}\mu_n(b)\, P^{(q)}_s(\f a b,\xi)\in\Z[\xi] $$

\no
et

$$d_n^{A}\mu_n(b)\, P^{(q)}_0(\f a b,\xi)\in\Z[\xi].$$
\no
\textbf{D\'emonstration}

\no D\'emontrons d'abord le premier et le troisi\`eme point.
Supposons $j\in[0,n-1]$ (le cas $j=n$ se traite de mani\`ere similaire, en se limitant \`a $s\leq q$)

\no
D'apr\`es (\ref{rjs}) $$
r^{(q)}_{j,s}(x)=\f{(-1)^{A-s}}{(A-s)!}\left(\f{\td{d}}{\td{d}l}\right)^{A-s}
\left[n!^{A-1}\f{(-l-x)_{n+1}}{(-l)^A_{n}(n-l)^q}(j-l)^A\right]_{|l=j}.$$
\'Ecrivons $$n!^{A-1}\f{(-l-x)_{n+1}}{(-l)^A_{n}(n-l)^q}(j-l)^A=F(l)\,G(l)^{A-1}\,H(l),$$
$$\text{o\`u }F(l)=\f{(-l-x)_{n}}{(-l)_{n+1}}(j-l),\quad G(l)=\f{n!} {(-l)_{n+1}}(j-l)
\quad\textup{ et }\ \ \ H(l)=(-l+n)^{A-q}(n-l-x) \,.$$

\sm
\no
D\'ecomposons maintenant  $F(l)$ et $G(l)$ en \'el\'ements simples

$$F(l)=1+\Sum{\begin{subarray}{c}
m\neq j \\ 0\leq m\leq n
\end{subarray}}{}\f{(j-m)f_{m}}{m-l},\quad G(l)=\Sum{\begin{subarray}{c}
m\neq j \\ 0\leq m\leq n
\end{subarray}}{}\f{(j-m)g_{m}}{m-l},$$

\sm

\no
o\`u

\sm

$$f_{m}=\f{(-m-x)_{n}}{\Prod{\begin{subarray}{c}
h\neq m \\ 0\leq h\leq n
\end{subarray}}{}{(h-m)}}=(-1)^m\f{(-m-x)_n}{n!}\binom{n}{m}$$

\sm

\no
et

\sm

$$g_{m}=\f{n!}{\Prod{\begin{subarray}{c}
h\neq m \\ 0\leq h\leq n
\end{subarray}}{}{(h-m)}}=(-1)^m\binom{n}{m}.$$

Il est imm\'ediat que $g_m$ est un entier. D'autre part $n!\,f_m\in\Z[x]$, donc $p^{\partieentiere{\f{n}{p-1}}}f_m\in\Z_p[x]$. De plus
le lemme \ref{aritmu} implique que pour $x=\f a{b}$, $\mu_n(b)f_{m}$ est un entier. On note $D_\lambda=\f1{\lambda!}\pa{\f\d{\d l}}^ \lambda$, on a alors
pour tout entier $\lambda\geq0$:

$$
(D_\lambda F(l))_{|l=j} =
\delta_{0,\lambda}-\displaystyle{\sum_{\begin{subarray}{c} m\,\neq j \ \ 0\leq
m\leq n
\end{subarray}}}\f{f_{m}}{(m-j)^{\lambda}} $$

\no
et

$$ (D_\lambda G(l))_{|l=j} =-\displaystyle{ \sum_{\begin{subarray}{c} m\neq j \ \ 0\leq m\,\leq n
\end{subarray}}{}\f{g_{m}}{(m-j)^{\lambda}}}.$$

\sm
\no
On a donc montr\'e que, pour tout $\lambda$ entier positif,

\begin{eqnarray}\label{GF}
\text{$d_n^\lambda\pa{D_\lambda G(l)}_{|l=j}\in\Z$ \ \ \ et \ \ \ $p^{\partieentiere{\f{n}{p-1}}}d_n^\lambda\pa{D_\lambda F(l)}_{|l=j}\in\Z_p[x].$}
\end{eqnarray}

\no
De plus, pour $x=\f a{b}$, $d_n^\lambda\mu_n(b)\pa{D_\lambda F(l)}_{|l=j}$ est entier. Enfin les d\'eriv\'ees $D_\lambda(H(l))_{|l=j}$ sont des polyn\^omes de $\Z[x]$ de degr\'e au plus $1$, et pour $x=\f a{b}$, $b\,D_\lambda(H(l))_{|l=j}$ est entier.

\no
Gr\^ace \`a la formule de Leibniz, on a

\sm

$$
D_{A-s}\left[n!^{A-1}\f{(-l-x)_{n+1}}{(-l)^A_{n}(n-l)^q}(j-l)^A\right]_{|l=j}=\sum_\nu
\pa{D_{\nu_0}\pa{F}}_{|l=j}
\pa{D_{\nu_1}\pa{G}}_{|l=j}\cdots.\pa{D_{\nu_{A-1}}\pa{G}}_{|l=j}\pa{D_{\nu_{A}}\pa {H}}_{|l=j}$$

\sm

\no
(o\`u la somme est sur les multi-indices $\nu\in\N^{A+1}$ tels que
$\nu_0+\cdots+\nu_{A}=A-s$), on d\'eduit alors que $p^{\partieentiere{\f n{p-1}}}d_n^{A-s}\,r^{(q)}_{j,s}(x)$ appartient \`a $\Z_p[x]$ et que
$b\,d_n^{A-s}\mu_n(b)\,r^{(q)}_{j,s}(x)$ est un \'el\'ement de $\Z$. Le premier et le troisi\`eme point sont alors d\'emontr\'es.

Pour le deuxi\`eme et le quatri\`eme point, en utilisant les \'equations (\ref{rfrac1}) et (\ref{rfrac2}), il suffit de montrer que

$$\f {p^{\partieentiere{\f{n}{p-1}}}d_n^{A-1}}{(A-1)!}\pa{\f{d}{\d l}}^{A-1}\pac{n!^{A-1}\f{(-l-x)_{n+1}}{(-l)^A_{n}(-l+n)^q}(j-l)^A
\f
{1}{l-k+x}}_{|l=j}\in\Z_p[x]$$

\no et, pour $x=\f a{b}$,

$$\frac{d_n^{A}\mu_n(b)}{(A-1)!}\pa{\f{d}{\d l}}^{A-1}\pac{n!^{A-1}\f{(-l-x)_{n+1}}{(-l)^A_{n}(-l+n)^q}(j-l)^A
\f
{1}{l-k+x}}_{|l=j}\in\Z.$$

\no
\'Ecrivons

\begin{eqnarray}
n!^{A-1}\f{(-l-x)_{n+1}}{(-l)^A_{n}(n-l)^q}(j-l)^A\f {1}{l-k+x} & =&\ti F(l)\,G(l)^{A-1}\,\ti H(l)
\end{eqnarray}

\no
avec
$$\ti F(l)=\f{(-l-x)_{n+1}}{(-l)_{n+1}}(j-l)\f {1}{l-k+x},\quad G(l)= \f{n!}{(-l)_{n+1}}(j-l),
\quad et \quad \ti H(l)=(-l+n)^{A-q}$$

\sm

\no
Gr\^ace au r\'esultat (\ref{GF}) sur $G$ et comme $D_\lambda(\ti H(l))_{|l=j} $ est un entier, il suffit de montrer que $p^{\partieentiere{\f{n}{p-1}}}d_n^{\lambda}D_\lambda(\ti F(l))_{|l=j}$ appartient \`a $\Z_p[x]$ et que pour $x=\f a{b}$, $d_n^{\lambda+1}\mu_n(b)D_\lambda(\ti F(l))_{|l=j}$ est entier; or

$$\ti F(l)=-1+\Sum{\begin{subarray}{c}
m\neq j \\ 0\leq m\leq n
\end{subarray}}{}\f{(j-m)\ti f_{m}}{m-l}$$

\no
avec

$$\ti f_{m}=\f{(-m-x)_{n+1}}{(m-k+x)\Prod{\begin{subarray}{c}
h\neq m \\ 0\leq h\leq n
\end{subarray}}{}{(h-m)}}=(-1)^m\f{(-m-x)_{n+1}}{n!(m-k+x)}\binom{n} {m}.$$

\no On voit donc que $p^{\partieentiere{\f n{p-1}}}\ti f_m$ est dans $\Z_p[x]$, et que, d'apr\`es le lemme \ref{aritmu}, pour $x=\f a{b}$, $d_n\mu_n(b)\ti f_m$ est entier.
La formule $$
(D_\lambda\ti F(l))_{|l=j} =
-\delta_{0,\lambda}-\displaystyle{\sum_{\begin{subarray}{c} m\,\neq j \ \ 0\leq
m\leq n
\end{subarray}}}\f{\ti f_{m}}{(m-j)^{\lambda}} $$
permet alors de conclure comme ci-dessus, et les deuxi\`eme et quatri\`eme points sont \'etablis.

\sm

\begin{cor}
Si $x$ est un nombre rationnel $\f a b$, alors $b\,
d_n^{A}\mu_n(b)\,S^{(q)}_{n}(x,\xi)$ est une combinaison lin\'eaire \`a
coefficients dans $\Z[\xi]$ de $\pa{\phi_s(x,\xi^{-1})}_{s\in[1,A]}\textup { et de }1$.
\end{cor}

\no
\textbf{D\'emonstration}
On utilise le corollaire \ref{corcombSxxi} et la proposition \ref{coefentier}.

\sm

\section{ Propri\'et\'es asymptotiques des polyn\^omes $P_s^{(q)}(x,z)$}

\begin{prop}\label{majcoef} Si $x$ est un nombre complexe fix\'e,
alors

$$\limsup_{n\rightarrow+\infty}{|P^{(q)}_s(x, \xi)|^{\f1n}}\leq 2^ {A-1}. $$
\end{prop}

\no
\textbf{D\'emonstration} Puisque
$$|P^{(q)}_s(x,\xi)|\leq\sum_{j=0}^n\abs{r^{(q)}_{j,s}(x)}$$ 
\no il nous suffit de majorer $r_{j,s}^ {(q)}(x)$. Or on
a

\begin{eqnarray}
r_{j,s}^{(q)}(x)&=&\f1{2\pi
i}\Int{|z+j+x|=\f12}{}{R^{(q)}_n(z)(z+j+x)^{s-1}\td{d}z}\\[4mm]
& =&\f1{2\pi i}\Int{|z+j+x|=\f12}{}
n!^{A-1}\f{(z)_{n+1}}{(z+x)^A_{n}(x+z+n)^q}
(z+j+x)^{s-1}\td{d}z
\end{eqnarray}

\no
On en d\'eduit que
\begin{eqnarray}
|r_{j,s}^{(q)}(x)|  &\leq&2^{-s} \sup_{|z+j+x|=\f12}\left(n!^{A-1}\frac{|(z)_{n+1}|}{| (z+x)^A_{n}(x+z+n)^q|}\right)\label{majrjs}
.\end{eqnarray}
\sm
Soit $m$ un entier positif tel que $\abs{x}+\frac12\leq m$, on a, pour $z$ tel que $|z+j+x|=\f1{2}$,

\begin{eqnarray}
|(z)_{n+1}| & = &\prod_{k=0}^{n}{|z+k|} \nonumber\\
    & = &\prod_{k=0}^{n}{|z+j+x -j- x+k|}\nonumber\\
& \leq & \prod_{k=0}^{n}{\left(\frac12+|x|+|k-j|\right)}\nonumber\\
& \leq & \prod_{k=0}^{n}{\left(m+|k-j|\right)}\nonumber\\
|(z)_{n+1}|  &\leq & m(m+1)...(m+j)(m+1)...(m+n-j)\leq (m+j)!(m+n-j)!.
\label{maj1}
\end{eqnarray}
Maintenant minorons

\begin{eqnarray}
|(z+x)_{n}| & = &\prod_{k=0}^{n-1}{|z+k+x|}\nonumber\\
    & = &\prod_{k=0}^{n-1}{|z+j+x -j+k|}\nonumber\\
& \geq & \prod_{k=0}^{n-1}{\left|-\frac12+|k-j|\right|} \nonumber
\end{eqnarray}

\no
En minorant par $\left||k-j|-\f12\right|\geq|k-j|-1$ si $|k-j|>1$, et par $\left||k-j|-\f12\right|\geq\f12$ sinon, on obtient, dans tous les cas
\begin{eqnarray}\label{maj3}
|(z+x)_{n}|\geq\f1{8n^3}j!(n-j)!.
\end{eqnarray}

\no En utilisant (\ref{maj1}) et (\ref{maj3}), on en d\'eduit que
$$
\f{|(z)_{n+1}|}{|(z+x)_n|}\leq8n^3\prod_{k=1}^m(j+k)(n-j+k)\leq8n^3(n+m)^{2m}.
$$  
\no Enfin
\begin{equation}
|(z+n+x)^q|=|(z+x+j-j+n)|^q\geq\f1{2^q}.
\label{maj4}\end{equation}

\no
On d\'eduit en utilisant 
(\ref{maj1}), (\ref{maj3}) et (\ref{maj4})  dans (\ref{majrjs}) que

$$\abs{r^{(q)}_{j,s}(x)}\leq2^{-s+q+3A}(n+m)^{2m}n^{3A}\binom{n}{j}^{A-1}
$$

\sm

\no
Comme $\binom{n}{j}\leq2^n$, il en r\'esulte que

$$|P^{(q)}_s(x,\xi)|\leq2^{-s+q+3A}(n+m)^{2m+3A+1}2^{n(A-1)}.$$
\no
On conclut donc

$$\limsup_{n\rightarrow+\infty}\abs{P^{(q)}_{s}(x,\xi)}^{\f1n}\leq 2^ {A-1}.$$

\begin{cor}\label{majplacinf} Pour tout $s\in[0,A]$, on a
\begin{eqnarray}
\limsup_n\f1n\ln\abs{b\mu_n(b)d_n^{A} P^{(q)}_s(x,\xi)}\leq \ln b+\sum_{q|b}\f{\ln q}{q-1}+
A+(A-1)\ln2.
\end{eqnarray}
\end{cor}

\section{Ind\'ependance lin\'eaire des formes lin\'eaires}

On consid\`ere la matrice

\begin{eqnarray}\label{det1}
M_n(x,z)=\pa{P^{(q)}_s(x,z)}_{\begin{array}{c}
q\in[0,A]\\
s\in[0,A]
\end{array}}
\end{eqnarray}
et on note

$$\Omega_n(x,z)=\det M_n\,.$$

\begin{prop} \label{indptlin}
On a

\begin{eqnarray}
\Omega_n(x,z)=\gamma
z^{n+1}(z-1)^{(A-1)n-2}{x^A}
\end{eqnarray}

\no
o\`u $\gamma\in\Q^*$.
\end{prop}
La preuve de cette proposition r\'esultera des lemmes suivants.

\begin{lem} \label{ordx0}

Le polyn\^ome $\Omega_n(x,z)$ est divisible par $x^A$.

\end{lem}

\no
\textbf{D\'emonstration} En d\'erivant (\ref{rjs}), on a 

$$\f{\d}{\dx}r^{(q)}_{j,1}(x)=\left\{\begin{array}{cl}
\f{(-1)^{A-1}}{(A-1)!}\left(\f{\td{d}}{\td{d}l}\right)^{A-1}
\left[n!^{A-1}\f{(-l-x)_{n+1}}{(-l)^A_{n}(n-l)^q}(j-l)^A\,\Sum{k=0}{n} \f{-1}{k-l-x}\right]_{|l=j}&
\textup{ si $j\in[0,n-1]$}  \\[6mm]
\f{(-1)^{q-1}}{(q-1)!}\left(\f{\td{d}}{\td{d}l}\right)^{q-1}
\left[n!^{A-1}\f{(-l-x)_{n+1}}{(-l)^A_{n}}\,\Sum{k=0}{n}\f{-1}{k-l-x}
\right]_{|l=n}&\textup{ si $j=n$ et $q>0$}\\[6mm]
0& \textup{ si $j=n$ et  $q=0$}
\end{array}\right.$$

\no
En utilisant la formule de Leibniz, on obtient, pour $j\in[0,n-1]$,

$$\f{\d}{\dx}r^{(q)}_{j,1}(x)$$ $$=$$

$$\f{(-1)^{A}}{(A-1)!}\Sum{u=0}{A-1}\binom{A-1}{u}\left(\f{\td{d}} {\td{d}l}\right)^{A-1-u}
\left[n!^{A-1}\f{(-l-x)_{n+1}}{(-l)^A_{n}(n-l)^q}(j-l)^A\right]_{|l=j} \left(\f{\td{d}}{\td{d}l}\right)^{u}\left[\Sum{k=0}{n}\f1{k-l-x}\right]_{|l=j}$$

$$=$$

$$\Sum{u=0}{A-1}(-1)^{u+1}r^{(q)}_{j,u+1}(x)\Sum{k=0}{n}\f1{(k-j- x)^{u+1}}.$$

\no
De m\^eme, pour $j=n$, on a

$$\f{\d}{\dx}r^{(q)}_{n,1}(x)= \left\{\begin{array}{cl}
\Sum{u=0}{A-1}{(-1)^{u+1}}r^{(q)}_{n,u+1}(x)\Sum{k=0}{n}\f1{(k- n-x)^{u+1}}&\textup{ si $q>0$}\\[4mm]
0& \textup{ si $q=0$.}\\

\end{array}\right.$$

\no
On en d\'eduit que

\begin{eqnarray*}
\f{\d}{\dx}r^{(q)}_{j,1}(x) & =& \Sum{k=0}{n}\f1{(k-j-x)^
{A}}\Sum{u=0}{A-1}{(-1)^{u+1}}r^{(q)}_{j,u+1}(x){(k-j-x)^{A-u-1}} \\
                    & =&\Sum{k=0}{n}\f1{(k-j-x)^{A}}\Sum{u=0}{A-1}{(-1)^{A-u}}r^{(q)}_{j,A-u}(x){(k-j-x)^{u}}.
\end{eqnarray*}

\no
Les quantit\'es que l'on d\'erive \'etant des polyn\^omes, les d\'eriv\'ees sont aussi des polyn\^omes. On en d\'eduit que, pour tout $j\in [0,n]$, pour tout $k\in[0,n]$, le polyn\^ome $\Sum{u=0}{A-1}{(-1)^{A-u}}r^ {(q)}_{j,A-u}(x){(k-j-x)^{u}}$ est divisible par $(k-j-x)^{A}$. Cela  implique que  pour tout $v\in[0,A-1]$, $\Sum{u=0}{v}{(-1)^{A-u}}r^{(q)}_{j,A- u}(x){(k-j-x)^{u}}$ est divisible par $(k-j-x)^{v+1}$ et donc que
$$\frac1{(k-j-x)^{v}} \Sum{u=0}{v}
{(-1)^{A-u}}r^{(q)}_{j,A-u}(x)(k-j-x)^{u}$$
 est un polyn\^ome  qui s'annule en
$x=k-j$.

\no
Il en r\'esulte en prenant $k=j$ que pour tout $j\in[0,n]$ et $v\in [0,A-1]$

$$\Sum{u=0}{v}
\f1{x^{v-u}}r^{(q)}_{j,A-u}(x)$$ est un polyn\^ome qui  s'annule en
$x=0$.

\no
Par lin\'earit\'e, on obtient que, pour tout $v\in[0,A-1]$

$$\Sum{u=0}{v}
\f1{x^{v-u}}P^{(q)}_{A-u}(x,z)=\Sum{u=0}{v}
\f1{x^u}P^{(q)}_{A-v+u}(x,z)$$ est un polyn\^ome en $x $  s'annulant en
$x=0$.

\no
Or, par multilin\'earit\'e sur les colonnes du d\'eterminant  (\ref{det1}), on obtient~:

\begin{eqnarray}
\Omega_n(x,z)=\abs{\begin{array}{ccccc}
P^{(0)}_0(x,z)&\Sum{u=0}{A-1} \f1{x^{u}}P^{(0)}_{u+1}(x,z) & \Sum {u=0}{A-2}
\f1{x^u}P^{(0)}_{u+2}(x,z)&\cdots &
P^{(0)}_A(x,z) \\
\cdots& \cdots & \cdots & \cdots\\
P^{(A)}_0(x,z)&\Sum{u=0}{A-1} \f1{x^u}P^{(A)}_{u+1}(x,z)& \Sum {u=0}{A-2}
\f1{x^u}P^{(A)}_{u+2}(x,z)&\cdots &
P^{(A)}_A(x,z) \\
\end{array}}
\end{eqnarray}

\vspace{0.5cm}
\no Chaque colonne (except\'e la premi\`ere) est un polyn\^ome en $x$ admettant $x=0$ comme z\'ero, on obtient donc que $x=0$
est z\'ero d'ordre au moins $A$ de $\Omega_n(x,z)$.  On conclut que
$x^A$ divise $\Omega_n(x,z)$.

\hspace{3mm}

\begin{lem}\label{degz}
Le polyn\^ome $\Omega_n(x,z)$ est de degr\'e $An-1$ en $z$.
\end{lem}

\no
\textbf{D\'emonstration}
En ajoutant \`a la premi\`ere colonne du d\'eterminant (\ref{det1}) les colonnes suivantes multipli\'ees par
$\phi_s\left(x,\frac1z\right)$, on obtient:

\begin{eqnarray}\label{det2}\Omega_n(x,z)=\abs{\begin{array}{cccc}
S^{(0)}_{n}(x,z)&P_1^{(0)}(x,z) &\cdots &P_A^{(0)}(x,z) \\
\cdots& \cdots & \cdots & \cdots\\
S^{(A)}_n(x,z)&P_1^{(A)}(x,z) &\cdots &P_A^{(A)}(x,z) \\
\end{array}}.
\end{eqnarray}

\sm

\no
Les \'el\'ements de la premi\`ere colonne sont exactement de degr\'e $-1$ en $z$, car le premier terme de la s\'erie $S^{(q)}_{n}(x,z)$ \'etant nul, on peut faire la somme \`a partir de $k=1$, et on obtient ainsi une s\'erie formelle en $\f1z$, de degr\'e $-1$. Les autres colonnes sont de degr\'e au plus $n$ en $z$, gr\^ace \`a la proposition \ref{approxpade}. On en d\'eduit que
le d\'eterminant est de degr\'e au plus $An-1$ en $z$. La proposition \ref{approxpade} nous montre
que les \'el\'ements surdiagonaux sont de degr\'e inf\'erieur ou \'egal \`a $n-1$ en $z$.
On en d\'eduit que dans le d\'eveloppement du d\'eterminant, tous les termes, autres que le produit des \'el\'ements diagonaux, sont de degr\'e  en $z$ strictement inf\'erieur \`a $An-1$. Mais les \'equations (\ref{rqq}) et (\ref{Pqs}) impliquent que $P^{(q)}_q(x,z)$ est exactement de degr\'e $n$ en $z$. Le produit des \'el\'ements diagonaux donne donc un \'el\'ement de degr\'e exactement $An-1$.  Le degr\'e en $z$ de
$\Omega_n(x,z)$ est donc exactement $An-1$.

\begin{lem}\label{majdegx}
Le polyn\^ome $\Omega_n(x,z)$ est de degr\'e au plus $A$ en $x$.
\end{lem}

\no
\textbf{D\'emonstration}
D\'eveloppons l'expression (\ref{det2}) du d\'eterminant $\Omega_n (x,z)$ par rapport \`a la premi\`ere colonne, on
obtient

$$\Omega_n(x,z)=\Sum{q=0}{A}{(-1)^q S^{(q)}_{n}(x,z)}\Delta_{q,0}(x,z)$$
o\`u les $\Delta_{q,0}(x,z)$ sont les d\'eterminants extraits.

\no
On a

\begin{eqnarray}\label{majdelta}
{x^{-A}\,S^{(q)}_{n}(x,z)}\Delta_{q,0}(x,z)=\Sum{k=0}{+\infty}x^{-A} \Delta_{q,0}(x,z)R^{(q)}_n(k)z^{-k}.\end{eqnarray}

\no
Cela implique, pour $\Re(x)>0$

\begin{eqnarray}
\abs{x^{-A}\Delta_{q,0}(x,z)R^{(q)}_n(k)}=\abs{x^{-A}\Delta_{q,0} (x,z)\,n!^{A-1}\frac{(k)_{n+1}}{(x+k)_n^A(x+k+n)^q}}
\leq n!^{A-1}\abs{\frac{(k)_{n+1}}{(x+n)^q}\frac{\Delta_{q,0}(x,z)}{x^ {A(n+1)}}}
.\end{eqnarray}

\no
La proposition \ref{approxpade} permet de majorer le degr\'e en $x$ de $\Delta_{q,0}(x,z)$ par $A(n+1)$. Cela implique que pour $z$ fix\'e quelconque, avec $\abs{z}>1$, $\frac {\Delta_{q,0}(x,z)}{x^{A(n+1)}}$ est born\'ee pour $\Re(x)>1$. On a donc pour $\Re(x)>1$

\begin{eqnarray}
\abs{x^{-A}\Delta_{q,0}(x,z)R^{(q)}_n(k)z^{-k}}\leq\frac K{\abs{x} ^q}(k)_{n+1}\abs{z^{-k}}
\label {n10000}
\end{eqnarray}
o\`u $K=K(z)$ est une constante ind\'ependante $x$.

\no
   Le terme de droite de l'in\'egalit\'e pr\'ec\'edente \'etant le terme g\'en\'eral d'une s\'erie convergente pour $\abs{z}>1$, on en d\'eduit que les termes de l'\'equation (\ref{majdelta}) tendent vers $0$ quand
$\Re(x)$ tend vers $+\infty$ si $q>0$ et restent born\'es pour $q=0$. On en conclut que le degr\'e en $x$ de $\Omega_n(x,z)$ est au plus $A$.

\vspace{4mm}

\begin{cor}\label{premiermorceaudelambda}
Le polyn\^ome  $\Omega_n(x,z)$ est de la forme

$${x^A}Q(z),$$

\no
o\`u $Q(z)$ un polyn\^ome de degr\'e $An-1$.

\end{cor}

\no
\textbf{D\'emonstration}
Cela r\'esulte des lemme \ref{ordx0}, \ref{degz} et \ref {majdegx}.

\begin{lem}\label{lambdaz0}
Le polyn\^ome $\Omega_n(x,z)$ est divisible par $z^{n+1}$.
\end{lem}
\no
\textbf{D\'emonstration}

\no
Du corollaire \ref{premiermorceaudelambda}, on d\'eduit

$$Q(z)=x^{-A}\Omega_n(x,z)=\lim_{\Re(x)\rightarrow+\infty}x^{-A} \Omega_n(x,z),$$

\no
d'o\`u

$$Q(z)=\Sum{q=0}{A}{(-1)^q \lim_{\Re(x)\rightarrow+\infty}x^{-A} S^ {(q)}_{n}(x,z)}\Delta_{q,0}(x,z).$$

\no
Le r\'esultat (\ref{n10000}) nous permet de conclure que pour $\abs{z}>1$, on a

$$Q(z)=\lim_{\Re(x)\rightarrow+\infty}x^{-A} S^{(0)}_{n}(x,z)\Delta_ {0,0}(x,z).$$

\no
On a de plus

\begin{eqnarray}\label{Qz}
\lim_{\Re{(x)}\rightarrow+\infty}x^{An}S^{(0)}_{n}(x,z)=\lim_{\Re{(x)} \rightarrow+\infty}n!^{A-1}\sum_{k=0}^{+\infty}\f{(k)_{n+1}x^{An}}{(k +x)^A_{n}}z^{-k}=n!^{A-1}\sum_{k=0}^{+\infty}(k)_{n+1}z^{-k}.\end{eqnarray}

\no
Or pour $\abs{Z}<1$, on a $$\sum_{k=0}^{+\infty}(k)_{n+1}Z^{k}=Z\f{\d^{n+1}}{\d Z^{n+1}}\sum_{k=0}^{+\infty}Z^{k+n}=Z\f{\d^{n+1}}{\d Z^{n+1}}\f{Z^n}{1-Z}=Z\f{\d^{n+1}}{\d Z^{n+1}}\f1{1-Z}=(n+1)!\f{Z}{(1-Z)^{n+2}}\cdot$$
Donc pour $\abs{z}>1$, 
$$\lim_{\Re{(x)}\rightarrow+\infty}x^{An}S^{(0)}_{n}(x,z)=\f{n!^A\,(n+1)\,z^{n+1}}{(z-1)^{n+2}}\cdot$$

Le fait que $\Delta_{0,0}$ soit un polyn\^ome en $x$ et $z$ de degr\'e au plus $A(n+1)$ en $x$ permet d'obtenir que

\begin{eqnarray}
\lim_{\Re{(x)}\rightarrow+\infty}x^{-A(n+1)}\Delta_{0,0}(x,z)
\end{eqnarray}

\no
est un polyn\^ome $M(z)$.

\no
On a donc $Q(z)=M(z)\f{n!^A\, (n+1)\,z^{n+1}}{(z-1)^{n+2}}$ et il en r\'esulte que $z^{n+1}$ divise $Q(z)$.

\sm

\begin{lem}\label{lambdaz1}
Le polyn\^ome $\Omega_n(x,z)$ est divisible par $(z-1)^{(A-1)n-2}$
\end{lem}

\no
\textbf{D\'emonstration}
\no
Pour $z\in\C\setminus]-\infty,0]$, on pose $z^{-t}=e^{-t\log z} $ o\`u $\log z$ est la d\'etermination du logarithme de $z$ de partie imaginaire comprise entre $-\pi$ et $\pi$.
Consid\'erons    l'int\'egrale

$$J_n^{(q)}(z)=\f1{2\pi i}\int_{\abs{t+x}=n+1}R^{(q)}_n(t)\,z^{-t} \dt$$
qui d\'efinit une fonction holomorphe pour $z\in\C\setminus]-\infty,0]$.

\sm \textbf{La nullit\'e en 1}

\sm

\no
Par d\'erivation sous le signe somme, on obtient

$$\f{\d^k J_n^{(q)}}{\dz^k}(z)=\f{(-1)^k}{2\pi i}
\int_{\abs{t+x}=n+1}R^{(q)}_n(t)\,(t)_k\,z^{-t-k}\dt.$$

\no
On remarque que
$$\deg_t R^{(q)}_n(t)\,(t)_k=n+1+k-An-q.$$

\no
Or l'int\'egrale d'une fonction rationnelle de degr\'e inf\'erieur ou \'egal \`a $-2$ est
nulle sur un contour ferm\'e  contenant l'ensemble de ses p\^oles. Cela implique que, si $$k\leq (A-1)n+q-3,$$ on a
 
\begin{eqnarray}\label{Jn}
\f{\d^kJ_n^{(q)}}{\dz^k}(1)=0.
\end{eqnarray}

\sm

\textbf{Lien avec $\Omega_n(x,z)$}

\sm
\no
La formule des r\'esidus nous donne

$$J_n^{(q)}(z)=\Sum{j=0}n \td{Res}_{[t=-j-x]}\left(R^{(q)}_n(t)\,z^ {-t}\right). $$

\no
On a

$$e^{-t\,\log z}=e^{(x+j)\,\log z}\Sum{k=0}\infty
\f{(-1)^k\,(t+x+j)^k(\log z)^k}{k!}.$$
\no
En utilisant les m\^emes notations que pour la proposition \ref {approxpade}, on obtient

$$\td{Res}_{[t=-j-x]}\left(R^{(q)}_n(t)\,z^{-t}\right)=
\Sum{s=1}{A}r^{(q)}_{j,s}(x)\f{(-1)^{s-1}\,e^{(x+j)\log z}(\log z)^ {s-1}}{{(s-1)}!}
.$$

\no
On en d\'eduit

\begin{eqnarray}
J_n^{(q)}(z)&=&\Sum{j=0}n\Sum{s=1}{A}r^{(q)}_{j,s}(x)\f{(-1)^{s-1}\,
e^{(x+j)\log z}\,(\log z)^{s-1}}{{(s-1)}!} \\
J_n^{(q)}(z)&=& e^{x\log z}\Sum{s=1}A P^{(q)}_s(x,z)\f{(-1)^{s-1}(\log z)^ {s-1}}{{(s-1)}!}
.\end{eqnarray}
\no
Dans (\ref{det2}), en ajoutant \`a la deuxi\`eme colonne les suivantes multipli\'ees respectivement par
$\f{(-1)^{s-1}(\log z)^{s-1}}{{(s-1)}!}$, on obtient

$$\Omega_n(x,z)=\abs{\begin{array}{ccccc}
S^{(0)}_n(x,z)& e^{-x\log z}J_n^{(0)}(z)  & P_2^{(0)}(x,z) &\cdots &P_A^{(0)}(x,z) \\
\cdots& \cdots & \cdots & \cdots & \cdots\\
S^{(q)}_n(x,z)& e^{-x\log z}J_n^{(A)}(z) & P_2^{(A)}(x,z) &\cdots &P_A^{(A)}(x,z) \\
\end{array}}$$
\no
Gr\^ace \`a (\ref{Jn}), les fonctions $J^{(q)}_n(z)$ ont un z\'ero en $z=1$  d'ordre au moins $(A-1)n-2$, cela
nous permet de conclure que $(z-1)^{(A-1)n-2}$ divise $\Omega_n(x,z)$.

\sm

\no
\textbf{D\'emonstration de la proposition \ref{indptlin}}

\no
Les lemmes \ref{lambdaz0} et \ref{lambdaz1} et le corollaire \ref {premiermorceaudelambda} permettent de conclure.

\section{Passage  du cas complexe au cas $p$-adique et d\'emonstration du th\'eor\`eme}

Pour $s$ complexe tel que $\Re(s)>1$, et $x$ r\'eel positif, on pose
$$T(s,x)=\f1{e^s}\sum_{l=0}^{e-1}\xi^{-l}\zeta(s,\frac{ x+l}  e).$$
\no Comme la fonction $s\longmapsto\zeta(s,\f{x+l}{e})$ peut \^etre prolong\'ee en une fonction
holomorphe sur  $\C\setminus\{1\}$ admettant le point $1$ p\^ole simple d'ordre $1$ et de r\'esidu
$1$, la fonction $s\longmapsto T(s,x)$ peut \^etre consid\'er\'ee comme
une fonction holomorphe sur $\C$.

Pour un nombre $p$-adique $x$, tel que $|x|_p\geq p$ et $s$ un entier strictement positif, on pose

\begin{eqnarray}
T_p(s,x)&=\displaystyle{\sum_{j=0}^{e-1}\f{\pa{\f {x+j} e}^{1-s}}
{e^s\left\langle \frac {x+j} e\right\rangle^{1-s}}\xi^{-j}{\zeta_p(s, \frac {x+j} e)} }&\textup{ si }s>1
\end{eqnarray}
et
\begin{eqnarray}
T_p(1,x)&= \displaystyle{\frac1e\lim_{s\rightarrow1}\sum_{j=0}^{e-1} \xi^{-j}{\zeta_p(t,\f {x+j} e)} }.
\end{eqnarray}

\no
On remarque que

\begin{eqnarray}
T_p(s,x)=\f1{e^s}\omega\pa{\f xe}^{1-s} \tilde T_p(s,x).
\end{eqnarray}

\begin{prop}\label{propfin}

Soit $x=\f{a}b$ un rationnel, tel que $|x|_p\geq p$, et soit $A$ un entier sup\'erieur ou \'egal \`a 2.

\no
Alors la dimension $\tau$ de l' espace vectoriel engendr\'e sur $\Q (\xi)$ par la famille
$\paa{1,\left(T_p(s,x)\right)_{s\in[1,A]}}$ v\'erifie

\medskip

$$\tau\geq\frac{[\Q_p(\xi):\Q_p]}{\varphi(e)}\frac{A\ln\abs{x}_p}{\ln b+\Sum{q| b}{}\f{\ln q}{q-1}+
A+(A-1)\ln2}.$$

\end{prop}

\begin{lem}\label{gasympreel}
Soient $s$ r\'eel, $s>1$, et $x$ r\'eel, $x>0$. On a pour tout entier $k>0$

$$T(s,x)=\f1{e(s-1)}\sum_{l=0}^{e-1}\xi^{-l}(x+l)^{1-s}-
\sum_{j=1}^{k-1}\binom{-s}{j-1}e^{j-1}\f{B_j}j\sum_{l=0}^{e-1}\xi^{-l}(x+l)^{1-s-j}+O_{x\rightarrow+\infty} (x^{-s-k+1})$$

\no
et
$$T(1,x)=-\frac{1}{e}\sum_{l=0}^{e-1}\xi^{-l}\ln(1+\frac l x)+
\sum_{j=1}^{k-1}e^{j-1}\frac{(-1)^jB_j}j\sum_ {l=0}^{e-1}\xi^{-l} (x+l)^{-j}
+O_{x\rightarrow+\infty} (x^{-k}).$$
\end{lem}
\no
\textbf{D\'emonstration}
Le cas $s>1$ est une application directe de l'\'equation (\ref {zetaasympt1}). Le cas $s=1$ r\'esulte de l'\'equation (\ref{zetaasympt2}) car on peut \'ecrire 
$$T(s,x)=\f1{e^s}\sum_{l=0}^{e-1}\xi^{-l}\left(\zeta(s,\frac{ x+l}  e)-\f1{s-1}\right).$$

\begin{lem} \label{gasymppadic}
Soient $p$ un nombre premier, $s$ un entier plus grand que $1$ et $x$ un \'el\'ement
de $\Q_p$, tel que $|x|_p\geq p$. On a
\begin{eqnarray}
T_p(s,x)=\f1{e(s-1)}\sum_{l=0}^{e-1}\xi^{-l}(x+l)^{1-s}-
\sum_{j=1}^{\infty}\binom{-s}{j-1}e^{j-1}\f{B_j}j\sum_{l=0}^{e-1}\xi^{-l}(x+l)^{1-s-j}
\end{eqnarray}
et
\begin{eqnarray}
T_p(1,x)&=&\displaystyle{-\frac1{e}\sum_{l=1}^{e-1}\xi^{-l}\log_p(1+ \frac lx)+\sum_{j=1}^{+\infty}e^{j-1}\frac{(-1) ^jB_j}j\sum_{l=0}^{e-1}\xi^{-l}(x+l)^{-j}}.
\end{eqnarray}

\end{lem}

\no
\textbf{D\'emonstration}
Le cas $s>1$ est une cons\'equence directe de l'\'equation (\ref {zetapasympt1}). Pour le cas $s=1$, en utilisant  l'\'equation (\ref {zetapasympt2}), on obtient:

\begin{eqnarray}
\lim_{s\rightarrow 1}\sum_{l=0}^{e-1}{\xi^{-l}\zeta_p(s,\frac {x+l} e)}& = &\sum_{l=0}^{e-1} \xi^{-l}\left(-\ln \left<\frac {x+l}e\right>+\sum_ {j=1}^{+\infty}{(-1)^jB_j\over j}\left(\frac{x+l}e\right)^{-j}\right)
\end{eqnarray}

\no
On doit distinguer deux cas

\begin{itemize}
\item
Si $\abs{x}_p\geq q_p$, alors pour $l$ compris entre 1 et $e-1$, on a

$$\omega\left(\frac{ x}e\right)=\omega\left (\frac{ x+l}e\right).$$

\no
On conclut que
\begin{eqnarray*}
\sum_{l=0}^{e-1}\xi^{-l}\log_p\left\langle \frac{ x+l}e\right\rangle&=& \sum_{l=0}^{e-1}\xi^{-l}\log_p\left( \frac{\frac{ x+l}e}{\omega\left (\frac{ x}e\right)}\right)\\
& = &\sum_{l=0}^{e-1}\xi^{-l}\log_p\left(1+\frac l x\right)+\sum_{l=0}^ {e-1}\xi^{-l}\log_p\left( \frac{ x}{e\,\omega\left(\frac{ x}e\right)} \right)\\
& = &\sum_{l=1}^{e-1}\xi^{-l}\log_p\left(1+\frac l x\right).
\end{eqnarray*}

\item Si $p=2$ et $\abs{x}_2=2$, on a alors 
$$\omega\pa{\frac{x+l}e}=(-1)^l\omega\pa{\frac{x}e}.$$

On conclut de m\^eme que
\begin{eqnarray*}
\sum_{l=0}^{e-1}\xi^{-l}\log_2\left\langle \frac{ x+l}e\right\rangle&=& \sum_{l=0}^{e-1}\xi^{-l}\log_2\left( \frac{\frac{ x+l}e}{(-1)^l\omega\left (\frac{ x}e\right)}\right)\\
& = &\sum_{l=0}^{e-1}\xi^{-l}\log_2\left(1+\frac l x\right)+\sum_{l=0}^ {e-1}\xi^{-l}\log_2\left( \frac{ x}{e\,\omega\left(\frac{ x}e\right)} \right)-\sum_{l=0}^{e-1}\xi^{-l}\log_2\left((-1)^l\right)\\
& = &\sum_{l=1}^{e-1}\xi^{-l}\log_2\left(1+\frac l x\right).
\end{eqnarray*}

\end{itemize}

\begin{prop} \label{majpadic}

Soit  $x=\f a b$ ($a$ et $b$ premier entre eux) un rationnel, tel que $|x|_p\ge p$, si on note

$$U^{(q)}_n(x)=b\mu_n(b)\,d_n^{A}\left(P^{(q)}_0(x,\xi)+\Sum{s=1}{A}P^ {(q)}_s(x,\xi)\,
T_p(s,x)\right),$$
\no
on a

$$\limsup_{n}\f1n\ln| U^{(q)}_n(x)|_p\leq  -A\ln|x|_p.$$

\end{prop}
\sm

\no
On va d\'emontrer cette proposition en plusieurs \'etapes.

\no
On pose

$$\tilde U^{(q)}_n(x)=d_n^{A}\left(P^{(q)}_0(x,\xi)+\Sum{s=1}{A}P^ {(q)}_s(x,\xi)\,
T_p(s,x)\right).$$

\begin{lem}\label{minordtildeu} On a, pour $x$ un nombre $p$-adique, tel que $|x|_p\geq p$
$$\tilde U^{(q)}_n(x)=\sum_{k=0}^{+\infty}\, u^n_k x^{-k}$$
\no
o\`u $(u^n_k)$  est une suite de nombres rationnels ind\'ependant de $x$, v\'erifiant  $$u^n_k=0$$ pour tout $k< A(n-1)-3$.

\end{lem}

\no
\textbf{D\'emonstration}
On a dans le corps des s\'eries de Laurent $\Q((1/x))$

$$(x+1)^{-j}=x^{-j}\sum_{m=0}^{+\infty}\binom {-j}{m} x^{-m}=\sum_{m=0}^{+\infty}\binom {-j}{m} x^{-m-j}$$

\no et

$$\log_p(1+\frac1x)=\sum_{m=1}^{+\infty}\frac{(-1)^{m+1}}m x^{-m},$$
On peut donc consid\'erer les s\'eries formelles de Laurent dans $\K((1/x))$ (o\`u $\K=\Q(\xi)$)
$$\Theta(s,x)=\f1{e(s-1)}\sum_{l=0}^{e-1}\xi^{-l}(x+l)^{1-s}-
\sum_{j=1}^{\infty}\binom{-s}{j-1}e^{j-1}\f{B_j}j\sum_{l=0}^{e-1}\xi^{-l}(x+l)^{1-s-j}$$
\no pour $s$ entier, $s>1$, et
$$\Theta(1,x)=\displaystyle{-\frac1{e}\sum_{l=1}^{e-1}\xi^{-l}\log_p(1+ \frac lx)+\sum_{j=1}^{+\infty}e^{j-1}\frac{(-1) ^jB_j}j\sum_{l=0}^{e-1}\xi^{-l}(x+l)^{-j}}.$$  

\no On peut calculer le terme g\'en\'eral de ces s\'eries en \'ecrivant, pour $s>1$,
$$
\Theta(s,x)=\f1e
\sum_{m=1}^{+\infty}\f1m\binom{-s}{m-1}\sum_{l=0}^{e-1}\xi^{-l}l^m x^{1-s-m}-
\sum_{j=1}^{\infty}\binom{-s}{j-1}e^{j-1}\f{B_j}j\sum_{l=0}^{e-1}\xi^{-l}\sum_{m=0}^{+\infty}\binom{1-s-j}ml^mx^{1-s-j-m}
$$ 
\no donc $$\Theta(s,x)=\sum_{k=0}^{+\infty}{a_{k,s}x^{-k}},$$ avec
$a_{k,s}=0$ pour $0\le k<s$ et, pour $k\ge s$~:
\begin{eqnarray}\label{gpssup1}
a_{k,s}=\f1{e(k-s+1)}\binom{-s}{k-s}\sum_{l=0}^{e-1}\xi^{-l}l^{k-s+1}-\sum_{j=1}^{k-s+1}e^{j-1}\f{B_j}j \binom{-s}{j-1}
\binom{1-s-j}{k-j-s+1}\sum_{l=0}^{e-1}\xi^{-l}l^{k-j-s+1}.
\end{eqnarray}
\no De m\^eme
$$\Theta(1,x)=\sum_{k=1}^{+\infty}a_{k,1}x^k$$ avec
\begin{eqnarray}\label{gps0}
a_{k,1}=\frac{(-1)^{k-1}}{ek}\sum_{l=1} ^{e-1}\xi^{-l} l^k+\sum_{j=1}^{k}(-1)^je^{j-1}\frac{B_j}j\binom{-j}{k-j}\sum_{l=0} ^{e-1}\xi^{-l} l^{k-j}.
\end{eqnarray}

\no En utilisant les lemmes \ref{gasympreel} et \ref{gasymppadic}, on voit que pour $x\in\Q_p$ avec $\abs{x}\ge p$, la s\'erie $\Theta(s,x)$ converge, et a pour somme $T_p(s,x)$, alors que pour $x$ r\'eel positif, on a pour tout entier $K\ge1$ le d\'eveloppement limit\'e
$$T(s,x)=\sum_{k=1}^{K-1}a_{k,s}x^{-k}+O(x^{-K}).$$
\no
Comme les polyn\^omes $P_s(x,z)$ sont de degr\'e au plus $n+1$ en $x $, cela implique que si on consid\`ere la s\'erie formelle dans $\K((1/x))$
$$V_n^{(q)}(x)=d_n^A(P_0^{(q)}(x,\xi)+\sum_{s=1}^AP_s^{(q)}(x,\xi)\Theta(s,x))=\sum_{k=-n}^{+\infty}\, u^n_k x^{-k}$$
on a pour $x$ r\'eel positif

$$d_n^AS_n^{(q)}(x,\xi))=\sum_{k=-n}^{K-1}\, u^n_k x^{-k}+O(x^{-K})$$
alors qu'au sens $p$-adique, pour $x$ rationnel tel que $\abs{x}\ge q_p$, $$\tilde U_n^{(q)}(x)=\sum_{k=-n}^{+\infty}\, u^n_k x^{-k}.$$
le corollaire \ref{corcombSxxi}  nous donne que

\begin{eqnarray}\tilde U^{(q)}_{n}(x)=o_{\Re(x)\to+\infty}(x^{-An+n+3-q}).\label{n100}
\end{eqnarray}

\no
L'unicit\'e du d\'eveloppement limit\'e montre donc que $u^n_k=0$ si $k< A(n-1)-3$.

\sm

\begin{lem}\label{majcoeftildeu} Les termes $u^n_k$ v\'erifient

$$\abs{u^n_k}_p\leq \frac{k+n+1}{\abs {e}_p}p^{\partieentiere{\frac{n}{p-1}}+1} .
$$
\end{lem}

\no
\textbf{D\'emonstration}
On rappelle que la valuation $p$-adique d'une nombre de Bernoulli est sup\'erieure ou \'egale \`a $-1$, par le th\'eor\`eme de Clausen-von Staudt (cf. \cite{Coh} pour une d\'emonstration) et que pour tout entier $n$ strictement positif et tout entier positif $i$, $\binom{-n}{i}=(-1)^i\binom{n+i-1}{i}$ est un entier. Les expressions (\ref{gpssup1}) et (\ref{gps0}) nous donnent donc directement, pour $s$ entier, $s\ge1$,
$$\abs{a_{k,s}}_p\leq\frac k{\abs {e}_p}p.$$
\no
La proposition \ref{approxpade} et la proposition \ref{coefentier}  assurent que les polyn\^omes $d_n^AP^{(q)}_s(x,z)$ sont de degr\'e  en $x$
au plus $n+1$ et ont des coefficients major\'es par $p^{\partieentiere {\frac{n}{p-1}}}$ en valeur absolue $p$-adique.

\no
On en d\'eduit, en utilisant la s\'erie formelle
$$V_n^{(q)}(x)=d_n^A(P_0^{(q)}(x,\xi)+\sum_{s=1}^AP_s^{(q)}(x,\xi)\Theta(s,x))=\sum_{k=-n}^{+\infty}\, u^n_k x^{-k}$$
que $$\abs{u^n_k}_p
\leq\frac{k+n +1}{\abs {e}_p}p^{\partieentiere{\frac{n}{p-1}}+1}.$$
\no

\begin{lem} \label{majpadictildeu}On a

$$\limsup_{n\rightarrow +\infty}\frac1n\ln\abs{\tilde  U^{(q)}_n(x)}_p\leq \frac{\ln p}{p-1}-(A-1)\ln\abs{x}_p$$
\end{lem}

\no
\textbf{D\'emonstration}
En utilisant le lemme \ref{minordtildeu}, on a

$$\abs{\tilde U_n^{(q)}(x)}_p\leq\sup_{k\geq(A-1)n-3}\abs{u^n_k}_p\abs {x}_p^{-k}.$$

\no
En utilisant le lemme \ref{majcoeftildeu}, on en d\'eduit que

$$\abs{\tilde U_n^{(q)}(x)}_p\leq \sup_{k\geq(A-1)n-3}\frac{k+n +1}{\abs {e}_p}p^{\partieentiere{\frac{n}{p-1}}+1}\abs{x}_p^{-k}$$

\no
On en d\'eduit que

$$\abs{\tilde U^{(q)}_n(x)}_p\leq M\sup_{k\geq(A-1)n-3}\pa{k+n +1}p^{\frac{n}{p-1}+1}\abs{x}_p^{-k}$$

\no
avec $M$ une constante ind\'ependante de $n$. La d\'ecroissance du terme de droite nous permet  alors de conclure pour $n$ suffisamment grand
$$\abs{\tilde U^{(q)}_n(x)}_p\leq M\pa{An-2}p^{\frac{n}{p-1}+1}\abs{x}_p^{-((A-1)n-3)}.$$

\vspace{4mm}
\no
\textbf{D\'emonstration de la proposition \ref{majpadic}}

\no
Comme $x=\f a b$  et $|x|_p\ge q_p$, on a $|b|_p=|x|_p^{-1}$. Il en
r\'esulte

\begin{center}{$|\mu_n(b)|_p=|x|_p^{-n}\,p^{-\pac{\f{n}{p-1}}}$}
\end{center}

\no
On en d\'eduit

\begin{eqnarray}
\lim_n\f1n\ln\abs{\mu_n(b)}_p=-\ln\abs{x}_p-\f{\ln p}{p-1}\label {normpadicmu}
\end{eqnarray}

\no
En utilisant le lemme \ref{majpadictildeu} et l'\'egalit\'e (\ref {normpadicmu}), on conclut

$$\limsup_n \f1n\ln\abs{U_n^{(q)}(x)}_p\leq \left(\frac{\ln p}{p-1}-(A-1)\ln{\abs{x}_p}\right)+\left(-\ln\abs{x}_p-\frac{\ln p}{p-1}\right)=-A\ln{\abs{x}_p}$$

\textbf{D\'emonstration du proposition \ref{propfin}}

La proposition \ref{coefentier} prouve que les coefficients des combinaisons lin\'eaires $U^{(q)}_n(x)$ sont des \'el\'ements de $ \mathcal{O}(\K)$. La proposition \ref{indptlin} donne l'ind\'ependance des formes lin\'eaires.  Le corollaire \ref{majplacinf}  donne une majoration de la valeur absolue aux places infinies des coefficients, ce qui permet de prendre

$$c= \ln b+\sum_{q|b}\f{\ln q}{q-1}+
A+(A-1)\ln2$$

\no  La proposition \ref{majpadic} nous donne une majoration de la valeur absolue $p$-adique des formes lin\'eaires. On prend

  $$\rho=A\ln{\abs{x}_p}$$

\no
  On peut donc appliquer le crit\`ere d'ind\'ependance lin\'eaire et la proposition est d\'emontr\'ee.

\vspace{1cm}

\no
\textbf{D\'emonstration du th\'eor\`eme 1}

\no
Le th\'eor\`eme 1 repose sur la proposition \ref{propfin} et le fait que 
$[\Q(\chi):\Q(\xi)]=\f{\varphi(v)}{\varphi(e)}$.

\vspace{1cm}

\no
\textbf{D\'emonstration du th\'eor\`eme 2}

\no
Appliquons la proposition \ref{propfin} \`a $x=\f2p$ et $e=2$.
On remarque alors que pour $p>2$, on a

$$\f{\pa{\f {\f2p+j} 2}^{1-s}}
{2^s\left\langle \frac {\f2p+j} 2\right\rangle^{1-s}}=\f1{2^s}.$$

\no
Cela implique

$$T_p\pa{s,\f2p}=\f1{2^s}\pa{\zeta_p\pa{s, \f 1 p}-\zeta_p\pa{s, \f {p+2}{ 2p}}}.$$ 

\no
On a de plus

\begin{eqnarray}\label{limdim}
\lim_{p\rightarrow+\infty}\frac{A\ln\abs{\f2p}_p}{\ln p+\f{\ln p}{p-1}+
A+(A-1)\ln2}=A.
\end{eqnarray}

\no
On en d\'eduit que pour $p$ suffisamment grand, l'espace vectoriel engendr\'e par
 
$$\pa{1,\pa{T_p\pa{s,\f2p}}_{s\in[1,A]}}$$
 est de dimension au moins $A$. L'équation (\ref{limdim}) permet de calculer explicitement la borne.
  Le th\'eor\`eme 2 est donc d\'emontr\'e. 

\vspace{1cm}
Je te tiens \`a remercier particuli\`erement Tanguy Rivoal pour m'avoir fourni son article sur les approximants de Pad\'e de la fonction de Lerch \cite{Ri2} et Henri Cohen pour l'aide que m'a apport\'e son enseignement sur les fonctions $p$-adiques.


\begin{thebibliography}{2}
   \bibitem[Be]{Be} F. Beukers, {\it Irrationality of some p-adic L- values}, (preprint)
   \bibitem[Ca]{Ca} F. Calegari, {\it Irrationality of certain p- adic periods for small p}, Intern. Math. Research Notices 2005:20 (2005), 1235--1249
   \bibitem[Coh]{Coh} H. Cohen, {\it Number Theory: Analytic and Modern Tools}, Springer, 2007
   \bibitem[Ma]{Ma} R. Marcovecchio, {\it Linear independence of forms in polylogarithms}
   \bibitem[Ri1]{Ri1} K. Ball and T. Rivoal, {\it Irrationalit\'e d'une infinit\'e de valeurs de la fonction z\^eta aux entiers impairs}, Invent. Math. 146:1 (2001), 193--207
  \bibitem[Ri2]{Ri2} T. Rivoal, { \it Simultaneous polynomial approximations of the Lerch function}, (preprint)


\end{thebibliography}
\end{document}